\documentclass{robincs}
\usepackage{prooftree,rotating,amsbsy}
\makeatletter

\let\orig@end@rotfloat = \end@rotfloat
\def\end@rotfloat{\orig@end@rotfloat\ignorespacesafterend}
\let\endsidewaysfigure\end@rotfloat
\let\endsidewaystable\end@rotfloat

\let\defn = \emph

\newcommand\alfa[3]{\alpha_{#1,#2,#3}}

\let\origint = \int
\def\int{{\textstyle\origint}}

\newcommand\curry{\mathrm{curry}}

\newcommand\Lin[1][\C]{\mathrm{Lin}_{#1}}

\@addtoreset{propn}{section}

\newcommand\pref[1]{\textup(\ref{#1}\textup)}

\newcommand\dlabel[2][]{\def\arg{#1}%
  \ifx\arg\@empty
    \refstepcounter{equation}%
	\global\let\robin@dnum = \theequation
  \else
    \def\@currentlabel{#1}%
	\gdef\robin@dnum{#1}%
  \fi
  \label{#2}}
\newcommand\dnum{\hbox\bgroup\morednum}
\newcommand\morednum[1][]{\def\arg{#1}\textup(%
	\ifx\arg\@empty%
		\robin@dnum%
	\else
		\ref{#1}%
	\fi
	\textup)\egroup}

\def\spleft#1{\hbox to 2em{\hss$#1$}}
\def\spright#1{\hbox to 2em{$#1$\hss}}

\newenvironment{mspill}{%
	\vskip\abovedisplayskip\hbox to \textwidth\bgroup\hss%
	\hbox\bgroup$\displaystyle
}{%
	$\egroup\hss\egroup\vskip\belowdisplayskip\noindent\ignorespacesafterend
}

\renewcommand\paragraph{\@startsection{paragraph}{4}{\z@}%
		{1ex}{-.5em}{\normalfont\normalsize\bfseries}}

\makeatother

\newif\ifincomplete
\incompletefalse

\newcommand\staraut{star-aut\-on\-om\-ous\xspace}
\newcommand\sm    {semi-mon\-oid\-al\xspace}
\newcommand\smc   {semi-mon\-oid\-al category\xspace}
\newcommand\smcs  {semi-mon\-oid\-al categories\xspace}
\newcommand\ssm   {symmetric semi-mon\-oid\-al\xspace}
\newcommand\ssmc  {symmetric semi-mon\-oid\-al category\xspace}
\newcommand\ssmcs {symmetric semi-mon\-oid\-al categories\xspace}
\newcommand\Ssmcs {Symmetric Semi-mon\-oid\-al Categories\xspace}
\newcommand\ssmcc {symmetric semi-mon\-oid\-al closed category\xspace}

\newcommand\Ssmccs{Symmetric Semi-mon\-oid\-al Closed Categories\xspace}

\newcommand\stac  {star-aut\-on\-om\-ous category\xspace}
\newcommand\sstac {semi star-aut\-on\-om\-ous category\xspace}
\newcommand\sstacs{semi star-aut\-on\-om\-ous categories\xspace}

\newcommand\SEC  {\textsc{ssmc}\xspace}

\newcommand\SECC {\textsc{ssmcc}\xspace}
\newcommand\SECCs{\textsc{ssmcc}s\xspace}

\newif\ifcomments\commentsfalse
\newcommand\comment[2][]{\ifcomments\marginpar{#1\raggedright\tiny#2}\fi}
\newcommand\dcomment{\comment[D]}
\newcommand\mll  {MLL$^-$\xspace}
\renewcommand\perp{^\bot}
\newcommand\perpp{\perp{}\perp}
\newcommand\dual[1]{#1\perp}

\newlength{\tw}
\newenvironment{plainbox}{\begin{center}\begin{tabular}{@{}|c@{\hspace{1ex}}|@{}}\hline
 \begin{minipage}[b]{\tw}\vspace*{.1in}
 \begingroup\small\baselineskip13pt}{\endgroup\end{minipage} \\[1ex] \hline
 \end{tabular}\end{center}}

\author{Robin Houston, Dominic Hughes and Andrea Schalk}
\title{Modelling Linear Logic without Units\\(Preliminary Results)}

\begin{document}
\maketitle

\section{Introduction}\label{intro}

Proof nets for \mll, unit-free Multiplicative Linear Logic
\citep{Girard87}, provide elegant, abstract representations of
proofs.
%
Cut-free \mll proof nets form a category, under path
composition\footnote{Path composition coincides with
normalisation by cut elimination \citep{Girard87}.}, in the manner of
Eilenberg-Kelly-Mac~Lane graphs \citep{EK66,KM71}.
Objects are formulas of \mll, and a morphism $A\to B$, a proof
net from $A$ to $B$, is a linking or matching between complementary
leaves (occurrences of variables).  For example, here is a morphism
from $p$ (a variable) to $p\otimes (q\tensor \dual{q})\perp$,
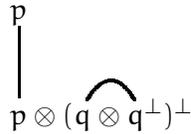
\[\thicklines
\begin{array}{l@{\hspace{10ex}}l}
        p\begin{picture}(0,0)\put(-3,-3){\line(0,-1){27}}\end{picture}
        \\[5ex]
        p\otimes(q%
        \begin{picture}(0,0)%
                \put(-2,9){\qbezier(0,0)(9,15)(18,0)}%
        \end{picture}%
        \otimes \dual{q})\perp
\end{array}\]
and here is an example of post-composing this proof net with another
from $p\otimes (q\tensor \dual{q})\perp$ to $\big((p\otimes
q)\perp\otimes q\big)\perp$, giving a proof net from $p$ to
$\big((p\otimes q)\perp\otimes q\big)\perp$:
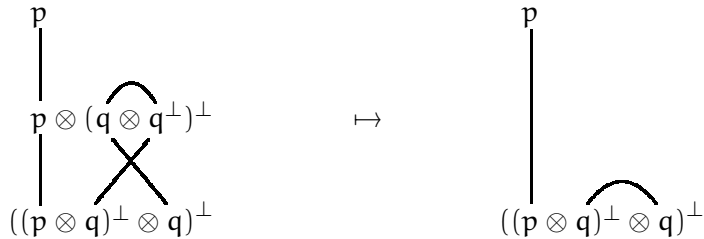
\[\thicklines
\begin{array}{l@{\hspace{10ex}}c@{\hspace{10ex}}l}
        p\begin{picture}(0,0)\put(-3,-3){\line(0,-1){27}}\end{picture}
        &&
        p\begin{picture}(0,0)\put(-3,-3){\line(0,-1){66}}\end{picture}
        \\[5ex]
        p\begin{picture}(0,0)\put(-3,-3){\line(0,-1){27}}\end{picture}
        \otimes(q
        \begin{picture}(0,0)%
                \put(-2,9){\qbezier(0,0)(9,15)(18,0)}%
                \put(0,-5){\qbezier(0,0)(10,-12)(20,-24)}%
        \end{picture}
        \otimes
        \begin{picture}(0,0)%
                \put(0,-5){\qbezier(0,0)(-10,-12)(-20,-24)}%
        \end{picture}
        \dual{q})\perp
        & \mapsto
        \\[4.6ex]
        \makebox[0pt][r]{$\big(($}p\otimes q)\perp\otimes q\big)\perp
        &&
        \makebox[0pt][r]{$\big(($}p\otimes q%
        \begin{picture}(0,0)%
                \put(-2,10){\qbezier(0,0)(13,17)(26,0)}%
        \end{picture}%
        )\perp\otimes q\big)\perp
\end{array}\]
To obtain the result of composition, one simply traces paths.

The category of proof nets is almost, but not quite, a \staraut
category \citep{BarrStac}.  The mismatch is the lack of units.
This
prompts the question of what categorical structure
\emph{does} match \mll.  More precisely:
\begin{quote}
\textbf{Question:}\\
	\it Can one axiomatise a categorical structure, a relaxation
	of \emph{star-autonomy}, suitable for modelling \mll?  The
	category of \mll proof nets \textup(with explicit negation\footnote{
	\emph{I.e.}, with formulas generated from variables by $\tn$
	and $(-)\perp$, rather than from literals by tensor and par.
	In other words, one drops the quotienting by de Morgan duality
	which is implicit in the usual definition of \mll formulas.}\textup)
	should be a free such category.
\end{quote}
Here by a \emph{categorical structure} we mean
\dcomment{Maybe you can improve the phrasing in this paragraph} 
data (functors, natural isomorphisms, and so forth) together with
coherence diagrams, akin to the axiomatisation of \staraut categories.

The naive proposal simply drops the units from a standard
axiomatisation of \staraut category $\C$ \citep{BarrStac}:
\begin{itemize}
\item Tensor. A functor $-\tensor-:\C\times \C\to\C$.
\item Associativity. 
A natural isomorphism 
\(
  \alpha_{A,B,C}:(A\tensor B)\tensor C\to A\tensor(B\tensor C)
\)
natural in objects $A,B,C\in\C$ such that the usual pentagon commutes.

\item \vspace*{-.5ex}Symmetry.
A natural isomorphism $\sigma_{A,B}:A\tn B\to B\tn A$ natural in
objects $A,B\in\C$ such that $\sigma_{B,A}^{-1}=\sigma_{A,B}$
and the usual hexagon commutes.

\item \vspace*{-1ex}Involution.
A functor $(-)\perp:\C\op\to\C$ together with a natural isomorphism
$A\to A\perpp$.
\item An isomorphism $\C(A\tensor B,C\perp)\to\C(A,(B\tn C)\perp)$
  natural in all objects $A,B,C$.
\end{itemize}
However, while there is a proof net from $p$ to $p\tn (q\tn\dual{q})\perp$
(the first proof net depicted at the beginning of the
Introduction), this axiomatisation fails to provide a corresponding
morphism from $p$ to $p\otimes (q\tensor \dual{q})\perp$ in the free
such category generated from the variables $p$ and $q$.
The problem of finding the right axiomatisation is
non-trivial.

\paragraph*{The solution predates the problem.} As so often, the solution 
long predates the problem.
\citet{DayPro} defines a \emph{promonoidal} category\footnote{
	The original paper on the subject \citep{DayPro}
	uses the term `premonoidal'. To add to the confusion,
	the word premonoidal is now used to mean something
	quite different.
} as a generalisation of a monoidal category.
Rather than having a functor \[\tn:\C\times\C\to\C\] and a unit object
$I\in\C$, a promonoidal category has functors
\[\begin{array}{c}
	P: \C\op\times\C\op\times\C\to\Set,\\
	J: \C \to \Set.
\end{array}\]
 This brings us to our primary definition:
\begin{plainbox}
\it
	A \textbf{\smc}
        $\C$ is a promonoidal category 
        such that\comment[D]{rephrased (a) for more precision, and (b) for typesetting prettiness}
	$P(A,B,C) = \C(A\tn B,C)$ for some functor $\tn:\C\times\C\to\C$.
\end{plainbox}
The motivation behind the choice of terminology \emph{semi} here is
twofold.  A monoidal category is a promonoidal category
satisfying:
\begin{itemize}
   \item[(a)] $P(A,B,C) = \C(A\tn B,C)$ for some functor $\tn:\C\times\C\to\C$, and
   \item[(b)] $J(A) = \C(I,A)$ for some object $I$.
\end{itemize}
\emph{Semi} refers to our use of just one of the two properties.
Secondly,
if one views \emph{semi} as short for \emph{semigroupal}, then one has
an analogy with semigroups, which are monoids without unit.

Emulating the usual progression from monoidal category to
\staraut category via symmetry and closure, we progress to
a notion of semi \staraut category, our candidate for an answer
to the question posed at the beginning of the Introduction.  (This
being a preliminary report on work in progress, we have yet to
complete the proof that the category of proof nets is a free semi
\staraut category.)
A key step in this progression is the following:
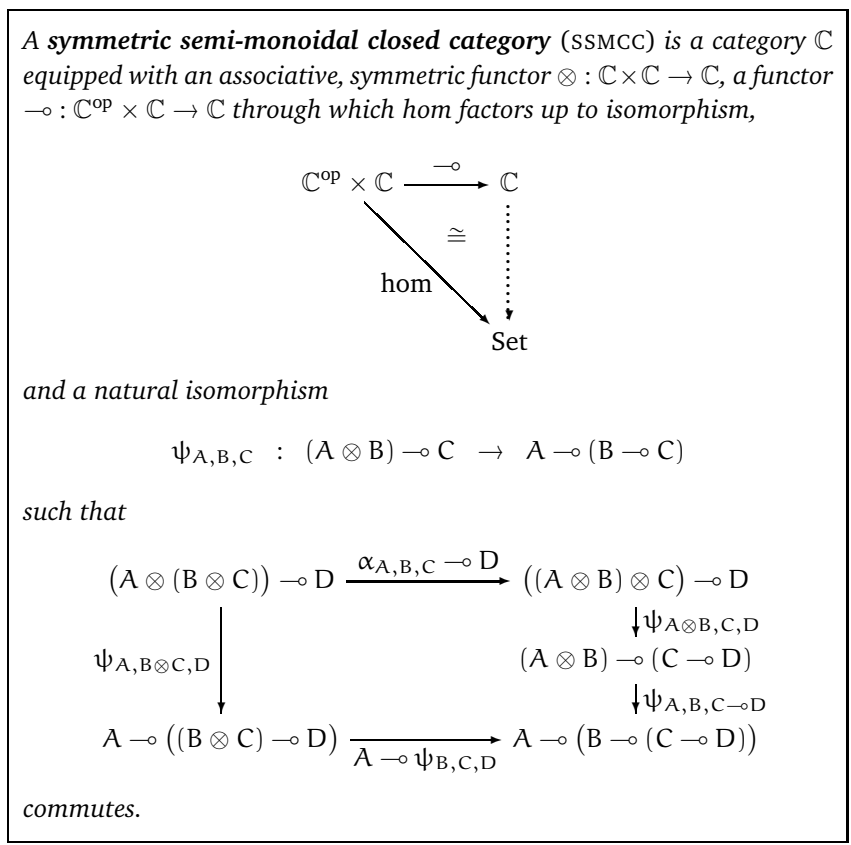
\begin{plainbox}
\it
	A \textbf{\ssmcc} \textup(\SECC\textup) is
	a category $\C$ equipped with an associative, symmetric
	functor
	$\tn:\C\times\C\to\C$, a functor $\lolli:\C\op\times\C\to\C$
        through which hom factors up to isomorphism,
        \begin{diagram}
        	\C\op\times\C &\rTo^\lolli& \C\\
        	&\rdTo_{\mathrm{hom}}\raise1em\hbox{\qquad$\cong$}&\dTo.\\
        	&&\Set
        \end{diagram}
        and a natural isomorphism
	\[
		\psi_{A,B,C}\;\;:\;\; (A\tn B)\lolli C \;\;\to\;\;
		A\lolli(B\lolli C)
	\]
	such that
	\begin{diagram}[h=1.5em] \bigl(A\tn(B\tn C)\bigr)\lolli D
	&\rTo^{\alpha_{A,B,C}\lolli D} &\bigl((A\tn B)\tn
	C\bigr)\lolli D\\ &&\dTo>{\psi_{A\tn B,C,D}}\\
	\dTo<{\psi_{A,B\tn C,D}}&&(A\tn B)\lolli(C\lolli D)\\
	&&\dTo>{\psi_{A,B,C\lolli D}}\\ A\lolli\bigl((B\tn C)\lolli
	D\bigr)&\rTo_{A\lolli\psi_{B,C,D}}
	&A\lolli\bigl(B\lolli(C\lolli D)\bigr) \end{diagram} commutes.
\end{plainbox}
See Section~\ref{s-sec} for details.  
Once again, 
the solution predates the problem: the natural isomorphism $\psi$, and
its commuting diagram are exactly as in the definition of symmetric
monoidal closed category in \citet{EKClosed}.

\paragraph*{Related work.}  
Our interest in obtaining an axiomatisation was sparked
by the desire to characterise the category of unit-free
proof nets for Multip\-lic\-at\-ive-Add\-it\-ive Linear Logic \citep{MallNets}
as a free category.

Soon after beginning to think about the problem, we came across an
interesting and informative proposal and discussion in a draft of
\citet{LSFreeBool}.
%
In this draft the authors define what they call \emph{(unitless)
autonomous categories}, motivated (like us) by the desire to model
unitless fragments of MLL.

Our definition of \SECC is apparently stronger, in that every \SECC is
a (unitless) autonomous category in the sense of (\emph{op.\ cit.}),
while the converse appears to be false.
In fact, certain properties are desired of categories in
\cite{LSFreeBool} which do not appear to be derivable from the given
conditions.  Indeed, in the presence of a tensor unit object $I$, the
axioms do not seem to imply symmetric monoidal closure.
One of the motivations behind producing this preprint is to suggest a
solution to the problem of finding a definition with the desired
properties.\footnote{In correspondence Lamarche and Stra\ss burger
have indicated that they may change the definition in the final
version of their paper, as a result.}
(Section~\ref{s-ls} discusses some other apparent divergences from the desired
properties.)

In our initial exploration of candidates for semi \staraut
category, we (independently) considered essentially the same
definition as in \cite{ProofNetCats}: a unitless linearly
distributive category with a suitable duality on objects.  Ultimately
we chose the approach more analogous to the standard progression from
monoidal category to \staraut category, via symmetry and
closure: it ties directly into the pioneering work of \citet{EKClosed},
with the \pref{psicoh} diagram, and also Day's promonoidal categories
\citep{DayPro}.

\paragraph*{Structure of paper.}
Section~\ref{s-summary} gives two different (but equivalent)
elementary definitions of \SECC.  Most of the remainder of the note is
devoted to showing that the first of these is equivalent to the conceptual
definition (Definition~\ref{def-secc}).

Section~\ref{s-bg} describes the background needed to understand the
subsequent development, in what is intended to be a clear and gentle
(but not rigorous) way. In particular, the definitions of \emph{coend}
and \emph{promonoidal category} are explained.
Section~\ref{s-sec} describes \ssmcs, using
Appendix~\ref{app-ax} (which gives a simple but
non-standard axiomatisation of symmetric monoidal categories).
Section~\ref{s-secc} shows that the conditions of \S\ref{s-summ-1}
are necessary and sufficient.
Section~\ref{s-sstac} defines \sstacs.

Finally section~\ref{s-related} discusses the relationship between
our definitions and the recent proposals of \citet{LSFreeBool}
and \citet{ProofNetCats}.

\paragraph*{Acknowledgements} 
Many thanks to Richard Garner and Martin Hyland for their penetrating
advice about an earlier version of this work. We are also grateful to
Lutz Stra\ss burger for useful discussion by email.  We acknowledge
use of Paul Taylor's diagrams package.


\section{Summary of Results}\label{s-summary}\dcomment{Perhaps we can skip this section, 
given the summary of the SSMCC definition in the Intro?  Entirely your call.}
The definition of \SECC (in terms of promonoidal categories) is
given in Def.~\ref{def-secc}. The main technical contribution of this
note is to show that this conceptual definition can be recast in
more elementary terms.
\subsection{First Description}\label{s-summ-1}
\begin{propn}\label{prop-1}
	A \ssmcc can be described by the following data:
	\begin{itemize}
	\item A category $\C$,
	\item Functors $\tn: \C\times\C\to\C$ and $\lolli:\C\op\times\C\to\C$,
	\item Natural isomorphisms $\alpha$, $\sigma$ and $\psi$ with
		components
		\[\begin{array}{r@{\;}l}
			\alpha_{A,B,C}:& (A\tn B)\tn C\to A\tn(B\tn C),\\
			\sigma_{A,B}:  & A\tn B\to B\tn A,\\
			\psi_{A,B,C}:  & (A\tn B)\lolli C \to A\lolli(B\lolli C),
		\end{array}\]
		such that
		\begin{starredequation}[$\boldsymbol\sigma$]\label{sigma}
		  \sigma_{A,B} = \sigma_{B,A}^{-1}
		\end{starredequation}
		  and the
		following coherence diagrams commute for all $A$,$B$,$C$,
		$D\in\C$:
		\begin{diagram}\dlabel[$\boldsymbol\alpha$]{pentagon}
		  \bigl((A\tensor B)\tensor C\bigl)\tensor D
		  &\rTo^{\alpha_{A\tn B,C,D}}&(A\tensor B)\tensor (C\tensor D)
		  &\rTo^{\alpha_{A,B,C\tn D}} & A\tensor \bigl(B\tensor (C\tensor D)\bigr)
		  \\
		  &\rdTo[snake=-1em](1,2)<{\alpha_{A,B,C}\tn D} &\dnum
		  & \ruTo[snake=1em](1,2)>{A\tn\alpha_{B,C,D}}
		  \\
		  & \spleft{\bigl(A\tensor(B\tensor C)\bigr)\tensor D}
		  & \rTo_{\alpha_{A,B\tn C,D}}
		  & \spright{A\tensor\big((B\tensor C)\tensor D\big),}
		\end{diagram}
		\begin{diagram}\dlabel[$\boldsymbol{\alpha\sigma}$]{hexagon}
		  (A\tn B)\tn C &\rTo^{\alfa ABC} & A\tn(B\tn C)
		    & \rTo^{\sigma_{A,B\tn C}} & (B\tn C)\tn A \\
		  \dTo<{\sigma_{A,B}\tn C} &&\dnum&& \dTo>{\alfa BCA}\\
		  (B\tn A)\tn C & \rTo_{\alpha_{B,A,C}} & B\tn(A\tn C)
		    & \rTo_{B\tn\sigma_{A,C}} & B\tn(C\tn A)
		\end{diagram}
		\begin{diagram}[h=1.5em]\dlabel[$\boldsymbol\psi$]{psicoh}
		\bigl(A\tn(B\tn C)\bigr)\lolli D &\rTo^{\alpha_{A,B,C}\lolli D}
			&\bigl((A\tn B)\tn C\bigr)\lolli D\\
		&&\dTo>{\psi_{A\tn B,C,D}}\\
		\dTo<{\psi_{A,B\tn C,D}}&\dnum&(A\tn B)\lolli(C\lolli D)\\
		&&\dTo>{\psi_{A,B,C\lolli D}}\\
		A\lolli\bigl((B\tn C)\lolli D\bigr)&\rTo_{A\lolli\psi_{B,C,D}}
			&A\lolli\bigl(B\lolli(C\lolli D)\bigr)
		\end{diagram}
	\item A functor $J:\C\to\Set$ and a natural isomorphism $e$
		with components $e_{A,B}: J(A\lolli B) \to \C(A,B)$.
	\end{itemize}
\end{propn}
This proposition is proved in section~\ref{s-secc} below.

\subsection{Second Description}\label{s-summ-2}
In fact there is a canonical choice for $J$ and $e$; in order to
state what it is, we need a few definitions:
\begin{definition}
  A \defn{category with tensor} is a category $\C$
  equipped with a functor $\tn:\C\times\C\to\C$ and
  a natural isomorphism $\alpha$ satisfying condition~\pref{pentagon}.

  A \defn{category with symmetric tensor} is a category $\C$
  with tensor, together with a natural isomorphism $\sigma$
  such that conditions \pref{sigma} and \pref{hexagon} hold.
\end{definition}
\begin{definition}\label{def-le}
	Let $\C$ be a category with tensor.
	A \defn{linear element} $a$ of the object $A\in\C$ is a natural transformation
	with components
	\[
		a_X: X \to A\tn X
	\]
	such that
	\begin{diagram}
	  X\tensor Y & \rTo^{a_X\tn Y} & (A\tn X)\tn Y\\
	  &\rdTo[snake=-1ex]_{a_{X\tn Y}} & \dTo>{\alpha_{A,X,Y}}\\
	  && A\tn (X\tn Y)
	\end{diagram}
	commutes for all $X$,$Y\in\C$.
\end{definition}
\begin{definition}\label{def-Lin}
	Given a category $\C$ with tensor, define a functor
	\[
		\Lin: \C\to\Set
	\]
	as follows. For $A\in\C$, let $\Lin(A)$ be the set of linear
	elements of $A$. For $f:A\to B$ and $a\in\Lin(A)$, let
	$\Lin(f)(a)$ be the linear element of $B$ with components
	\[
		X \rTo^{a_X} A\tn X \rTo^{f\tn X}B\tn X.
	\]
\end{definition}
It turns out that, in the situation of Prop.~\ref{prop-1}, there is a
canonical natural isomorphism between $J$ and $\Lin$. Furthermore, it
happens that this natural isomorphism takes $e$ to a particular natural
transformation $l$, which is defined as follows. 
\begin{definition}\label{def-l}
	Suppose we have $(\C,\tn,\alpha)$ as above, together with a functor
	\(
		\lolli: \C\op\times\C\to\C
	\)
	and a natural isomorphism
	\[
		\curry_{A,B,C}: \C(A\tn B,C) \to \C(A, B\lolli C)
	\]
	with counit (i.e.\ evaluation map) $\e^A_B:(A\lolli B)\tn A\to B$.
	
	Define the natural transformation $l_{A,B}: \Lin(A\lolli B) \to \C(A,B)$
	as follows: for each $x\in\Lin(A\lolli B)$, let $l_{A,B}(x)$ be the
	composite
	\[
		A \rTo^{x_{A}} (A\lolli B)\tn A  \rTo^{\e^A_B} B.
	\]
\end{definition}
Our first description (Prop.~\ref{prop-1}) is equivalent to the following.   
\begin{propn}\label{prop-2}
	An \SECC can be described by:
	\begin{itemize}
	\item A category $\C$ with symmetric tensor,
	\item A functor $\lolli:\C\op\times\C\to\C$ with a natural isomorphism
		\[
			\curry_{A,B,C}:\C(A\tn B,C) \to \C(A, B\lolli C)
		\]
	\end{itemize}
	such that the natural transformation
	$l$ of Def.~\ref{def-l}
	is invertible.
\end{propn}
The proof of Prop.~\ref{prop-2} is still in draft form,
and is not included in this preliminary note.

\section{Technical Background}\label{s-bg}
This section gives an informal introduction to coends and
promonoidal categories.
\subsection{Coends}\label{s-coends}
The definition of promonoidal category involves coends, so we
need to understand them to some extent. Fortunately they are
quite simple: a coend is just a slightly more general version
of a colimit.

Recall that a colimit of a functor $J: \D\to\C$ is a universal natural
transformation from $J$ to some object $X\in\C$. Coends just
extend this idea to mixed-variance functors $J: \D\times \D\op\to\C$:
a coend of $J$ is a universal \emph{dinatural}\footnote{
	The dinatural transformations we need to consider are
	of the special kind called \emph{extraordinary natural
	transformations}. The distinction is important in
	enriched category theory, where extraordinary natural
	transformations can be defined but dinatural transformations
	can not in general.
}
transformation from $J$ to an object $X\in\C$.

A dinatural transformation $\gamma: J\To X$ is a family
of arrows
\[
	\gamma_A: J(A,A) \to X
\]
indexed by the objects $A\in \D$, such that for every
$f:A\to B$ in $\D$ the diagram
\begin{diagram}
	J(A,B) & \rTo^{J(A,f)}& J(A,A)\\
	\dTo<{J(f,B)}&&\dTo>{\gamma_A}\\
	J(B,B) & \rTo_{\gamma_B}& X
\end{diagram}
commutes. Such a dinatural transformation is universal if, for every
object $Y$ and dinatural transformation $\delta: J\To Y$,
there is a unique fill-in morphism $g: X\to Y$ such that
\begin{diagram}[h=1.5em]
	&&X\\
	&\ruTo^{\gamma_A}\\
	J(A,A)&&\dTo>g\\
	&\rdTo_{\delta_A}\\
	&&Y
\end{diagram}
commutes for every $A\in\D$. 

It is conventional, and very handy, to write coends using
integral notation. The coend of $J$ is written $\int^{A\in\D} J(A,A)$;
we usually omit the `$\in\D$' part when it is obvious from the context.

In the rest of this note, we only need to use coends in $\Set$ or in
functor categories of the form $[\C,\Set]$ for some $\C$. These categories
have all (small) coends, so we shall never have to worry about whether
or not a particular coend exists.\footnote{
	We are glossing over some size issues here,
	which can be dealt with in the usual way.
	The foundationally conservative reader may
	read the word `category', where it appears
	in a definition, as `small category'.
}

Here are some important properties of coends.
We use them heavily in the sequel, often without remark.
\begin{itemize}
	\item Left adjoints preserve coends. In particular, in a cartesian
		closed category $\C$ we have
		\[
			A\times\int^{X\in\D} J(X,X) \cong \int^{X\in\D} A\times J(X,X)
		\]
		for every $A\in\C$;
	\item The `Fubini theorem': if $\int^{Y\in\D_2}J(X,X',Y,Y)$ exists for
			all $X$, $X'\in\D_1$ then
	\[
		\int^{X\in\D_1}\int^{Y\in\D_2} J(X,X,Y,Y) \cong
			\int^{(X,Y)\in\D_1\times\D_2}J(X,X,Y,Y)
	\]
	for $J: \D_1\times\D_1\op\times\D_2\times\D_2\op\to\C$;
	\item For every $F:\C\to\Set$ and $Y\in\C$,
	\[
		FY \cong \int^{X\in\C} FX\times\C(X,Y)
	\]
	and the coend on the right exists;
	\item The dual of the above: for every $F: \C\op\to\Set$ and $X\in\C$,
	\[
		FX \cong \int^{Y\in\C} FY\times\C(X,Y).
	\]
\end{itemize}
The latter two isomorphisms can be regarded as versions of the
Yoneda lemma. For proofs of all these facts, and a generally very
nice tutorial introduction to coends, see the lecture notes by
\citet{CatNotes}.
%

\subsection{Promonoidal Categories}\label{s-promon}
As mentioned in the introduction, a promonoidal category is a
category $\C$ together with functors
\[\begin{array}{c}
	P: \C\op\times\C\op\times\C\to\Set,\\
	J: \C \to \Set.
\end{array}\]
and natural isomorphisms $\alpha$, $\lambda$ and $\rho$
satisfying conditions analogous\footnote{
	In fact this is no mere analogy: monoidal and promonoidal
	categories are both instances of the general notion of
	\emph{pseudomonoid} in a monoidal bicategory \citep{MonBicat}.
	A promonoidal category is a pseudomonoid in the bicategory
	whose objects are categories and whose 1-cells are modules
	(aka profunctors or distributors).
} to those in the definition of a monoidal category.

In order to understand what these conditions are, we develop
an informal procedure for translating the language
of monoidal categories into the language of promonoidal categories.
In a monoidal category $\C$ we can form various functors
$\C^n\to\C$ for some natural number $n$, using the tensor
product and unit object. For example we have the functor
$F:\C^2\to\C$ defined as
\[
	F(A_1, A_2) = (A_1 \tn I) \tn (A_2 \tn A_1).
\]
In general, such a functor is described by an expression
formed from variables $A_1, \dots, A_n$ and the constant $I$
using the binary operation $\tn$. These expressions
can be thought of as denoting the objects of the free monoidal category
on a countably infinite set of generators.

\def\qq#1{\raise1pt\hbox{$\ulcorner$}\!#1\!\raise1pt\hbox{$\urcorner$}}%
Given such an expression $S$ and promonoidal category $\C$, we can
define a corresponding expression $\qq{S}$, representing a
functor $(\C\op)^n\to[\C,\Set]$, recursively as follows:
\[\begin{array}{r@{{}={}}l}
	\qq{I}     & J,\\
	\qq{A_i}   & \C(A_i, -),\\
	\qq{S\tn T}& \int^{X,Y} \qq S(X) \times \qq T(Y) \times P(X,Y,-).
\end{array}\]
(We assume that bound variables are renamed where necessary.)

For example, we have
\[\begin{array}{r@{\;}c@{\;}l}
	\qq{A_1\tn A_2} &=& \int^{X,Y} \qq{A_1}(X)\times\qq{A_2}(Y)\times P(X,Y,-)\\
		&=&     \int^{X,Y} \C(A_1,X)\times \C(A_2,Y)\times P(X,Y,-)\\
		&\cong& \int^{X} \bigl(\C(A_1,X)\times \int^Y \C(A_2,Y)\times P(X,Y,-)\bigr)\\
		&\cong& \int^{X} \C(A_1,X)\times P(X,A_2,-)\\
		&\cong& P(A_1,A_2,-)
\end{array}\]
and
\[\begin{array}{r@{\;}c@{\;}l}
	\qq{A_1\tn I} &=& \int^{X,Y} \qq{A_1}(X)\times\qq{I}(Y)\times P(X,Y,-)\\
		&=&     \int^{X,Y} \C(A_1,X)\times J(Y)\times P(X,Y,-)\\
		&\cong& \int^Y J(Y)\times P(A_1,Y,-).
\end{array}\]
Similarly we have
\[\begin{array}{r@{{}\protect\cong{}}l}
	\qq{I\tn A_1} & \int^X J(X)\times P(X,A_1,-),\\
	\qq{(A_1\tn A_2)\tn A_3} & \int^X P(A_1,A_2,X)\times P(X,A_3,-),\\
	\qq{A_1\tn (A_2\tn A_3)} & \int^X P(A_1,X,-)\times P(A_2,A_3,X).
\end{array}\]
Now we can say what the types of $\alpha$, $\lambda$ and $\rho$
should be. They should have components
\[\begin{array}{r@{\;}l}
	\alfa ABC:& \int^X P(A,B,X)\times P(X,C,-) \to \int^X P(A,X,-)\times P(B,C,X)\\
	\lambda_A:& \int^X J(X)\times P(X,A,-) \to \C(A,-)\\
	\rho_A:& \int^Y J(Y)\times P(A,Y,-) \to \C(A,-)
\end{array}\]
Each component here is \emph{itself} a natural transformation, so
we can add an extra variable and ask for natural isomorphisms
\[\begin{array}{r@{\;}l}
	\alpha_{A,B,C,Z}:& \int^X P(A,B,X)\times P(X,C,Z) \to \int^X P(A,X,Z)\times P(B,C,X)\\
	\lambda_{A,Z}:& \int^X J(X)\times P(X,A,Z) \to \C(A,Z)\\
	\rho_{A,Z}:& \int^Y J(Y)\times P(A,Y,Z) \to \C(A,Z)
\end{array}\]
between functors to $\Set$.

We impose the usual coherence conditions, as described for monoidal
categories in the appendix and elaborated below in the \sm
case.

\subsection{Notation for Diagrams}
In a diagram containing several cells which are known to commute,
we often label each such cell with the reason that it commutes,
removing the need for separate explanations that must be
cross-referenced with the diagram.
The symbol $\natural$ is used to indicate that a cell commutes
by naturality of some natural transformation.

\section{\Ssmcs}\label{s-sec}
We consider the special case of a \defn{\smc},
i.e.\ a promonoidal category $\C$ whose $P$ is represented by
a functor $\tn:\C\times\C\to\C$.
Therefore suppose we have such a functor $\tn$, and that
$P(A,B,X) = \C(A\tn B,X)$. Now we have
\[\begin{array}{r@{\;}c@{\;}l}
	\qq{(A_1\tn A_2)\tn A_3}  &\cong& \int^X \qq{A_1\tn A_2}(X)\times\C(X\tn A_3,-)\\
		&\cong& \int^X \C(A_1\tn A_2, X)\times\C(X\tn A_3,-)\\
		&\cong& \C((A_1\tn A_2)\tn A_3, -)
\end{array}\]
and similarly
\[
	\qq{A_1\tn(A_2\tn A_3)} \cong \C(A_1\tn (A_2\tn A_3), -).
\]
Thus, by Yoneda, the associativity isomorphism may be represented
by a natural isomorphism $\alpha$ with components
\[
	\alpha_{A,B,C}: (A\tn B)\tn C \to A\tn(B\tn C),
\]
just as in an ordinary monoidal category, subject to the
usual pentagon condition~\pref{pentagon}.

We are really interested in \emph{symmetric} \smcs, so suppose
that there is also a symmetry
$\sigma$ with components $\sigma_{A,B}: A\tn B \to B\tn A$
such that $\sigma_{A,B} = \sigma_{B,A}^{-1}$ for all $A$,$B\in\C$,
and satisfying the hexagon condition~\pref{hexagon}.

By the argument in the appendix -- more precisely, by reinterpreting
that argument in the promonoidal setting -- we have a \ssmc
just when there is a natural isomorphism
$\lambda$ with components
\[
	\lambda_{A,Z}: \int^X J(X)\times \C(X\tn A,Z) \to \C(A,Z)
\]
such that the diagram
\begin{diagram}\dlabel{bigone}
	\int^X J(X)\times\C((X\tn B)\tn C, Z) & \rTo^{\int^X J(X)\times\C(\alfa XBC, Z)}
		& \int^X J(X)\times\C(X\tn (B\tn C), Z)
	\\
	\dTo<{\cong}
	\\
	\int^{X,Y} J(X)\times\C(X\tn B, Y)\times \C(Y\tn C, Z)&\dnum
		&\dTo>{\lambda_{B\tn C,Z}}
	\\
	\dTo<{\int^Y\lambda_{B,Y}\times\C(Y\tn C,Z)}
	\\
	\int^Y\C(B,Y)\times\C(Y\tn C,Z) &\rTo_\cong&\C(B\tn C,Z)
\end{diagram}
commutes for all $B$,$C$,$Z\in\C$.

\section{\Ssmccs}\label{s-secc}
Generally speaking, a promonoidal category $\C$ is \defn{left closed}
if it has a functor $\lolli:\C\op\times\C\to\C$ with a natural
isomorphism
\[
    P(A,B,C) \cong \C(A, B\lolli C).
\]
In the case of present interest, we have:
\begin{definition}\label{def-secc}
	A \defn{\ssmcc} (\defn{\SECC})
	is a \ssmc $\C$ together with a functor
	\[
		\lolli:\C\op\times\C\to\C
	\]
	and a natural isomorphism with components
	\[
		\curry_{A,B,C}:\C(A\tn B, C)\to\C(A, B\lolli C).
	\]
\end{definition}
%
%
Recall the characterisation claimed in Prop.~\ref{prop-1}.
This differs from Def.~\ref{def-secc} in the following ways:
instead of $\lambda$ we have $e$;
instead of $\curry$ we have $\psi$; and
instead of \pref{bigone} we have \pref{psicoh}.
(Note that no conditions are imposed explicitly on $J$ or $e$.)
The rest of the section is devoted to proving Prop.~\ref{prop-1}.

\begin{lemma}\label{lemma-1}
	There is an isomorphism \SECC,
	$\qq{I\tn A} \cong J(A\lolli -)$, natural in $A$.
\end{lemma}
\begin{proof}
We have the chain of natural isomorphisms
\[\begin{array}{r@{\;}c@{\;}l}
	\qq{I\tn A} &=& \int^X J(X)\times \C(X\tn A,-)\\
		&\cong& \int^X J(X)\times \C(X,A\lolli -)\\
		&\cong& J(A\lolli -);
\end{array}\]
\vskip -\belowdisplayskip
\vskip -\baselinestretch\baselineskip
\end{proof}
\newcommand\JX{\mathchoice{\int^X JX\times}{\int^X JX\times}%
	{\int^X\!\!JX\times}{\int^X\!\! JX\times}}%
\newcommand\JXY{\mathchoice{\int^{X,Y} JX\times}{\int^{X,Y} JX\times}%
	{\int^{X,Y}\!JX\times}{\int^{X,Y}\! JX\times}}%
\begin{lemma}\label{lemma-2}
	To give a natural isomorphism $\lambda_A: \qq{I\tn A}\To\qq{A}$
	such that \pref{bigone} commutes
	is to give a natural isomorphism
	\[
		e_{A,Z}: J(A\lolli Z) \to \C(A, Z).
	\]
	such that the diagram
	\begin{diagram}[h=1.5em,labelstyle=\scriptstyle]\dlabel{diag-e}
		\JX\C(X\tn(B\tn C),Z) &\rTo^{\JX\C(\alpha_{X,B,C}, Z)} & \JX\C((X\tn B)\tn C,Z)\\
		&								& \dTo>{\JX\curry_{X\tn B,C,Z}}\\
		\dTo<{\JX\curry_{X,B\tn C,Z}}&	& \JX\C(X\tn B,C\lolli Z)\\
		 &								& \dTo>{\JX\curry_{X,B,C\lolli Z}}\\
		\JX\C(X,(B\tn C)\lolli Z)&\dnum	& \JX\C(X,B\lolli(C\lolli Z))\\
		\dTo<{\cong}&					& \dTo>{\cong}\\
		J((B\tn C)\lolli Z) &			& J(B\lolli(C\lolli Z))\\
		\dTo<{e_{B\tn C,Z}}&			& \dTo>{e_{B, C\lolli Z}}\\
		\C(B\tn C,Z)&\rTo_{\curry_{B,C,Z}}& \C(B,C\lolli Z)
	\end{diagram}
	commutes for all $B$,$C$,$Z\in\C$.
\end{lemma}
\begin{proof}
By Lemma~\ref{lemma-1} we can derive $e$ from $\lambda$ and
vice versa, and the translations are mutually inverse. So it
makes no difference whether we are given $e$ or $\lambda$,
since each can be derived from the other in a canonical way.

It remains to show that \pref{bigone} is equivalent to
\pref{diag-e}.
\begin{sidewaysfigure}
\[\begin{array}{@{\hskip4.5em}l}
\begin{diagram}[labelstyle=\scriptstyle]
	 \JX\C((X\tn B)\tn C,Z) &\rTo^{\JX\curry_{X\tn B,C,Z}}& \JX\C(X\tn B,C\lolli Z)
	\\
	 \dTo<\cong&&\dTo<\cong&\rdTo^{\JX\curry_{X,B,C\lolli Z}}
	\\
	 \rnode{a}{\JXY\C(X\tn B,Y)\times\C(Y\tn C,Z)}
		&\rTo_{\JXY\C(X\tn B,Y)\times\curry_{Y,C,Z}}
		& \JXY\C(X\tn B,Y)\times\C(Y,C\lolli Z)
		&\qquad&\JX\C(X,B\lolli(C\lolli Z))
	\\
	 \dTo>{\JXY\curry_{X,B,Y}\times\C(Y\tn C,Z)}&
	 	&\dTo<{\JXY\curry_{X,B,Y}\times\C(Y, C\lolli Z)}&\ldTo_\cong
	\\
	 \JXY\C(X,B\lolli Y)\times\C(Y\tn C,Z)
	 	&\rTo_{\JXY\C(X,B\lolli Y)\times\curry_{Y,C,Z}}
		&\JXY\C(X,B\lolli Y)\times\C(Y,C\lolli Z)
	\\
	 \dTo<\cong && \dTo<\cong && \dTo>\cong
	\\
	 \int^YJ(B\lolli Y)\times\C(Y\tn C,Z)
		&\rTo_{\int^Y\!\!J(B\lolli Y)\times\curry_{Y,C,Z}}
		&\int^YJ(B\lolli Y)\times\C(Y,C\lolli Z)
	\\
	 \dTo>{\int^Y\!\!e_{B,Y}\times\C(Y\tn C,Z)} &
		& \dTo<{\int^Y\!\!e_{B,Y}\times\C(Y,C\lolli Z)}
		&\luTo[leftshortfall=1.5em,crab+]^\cong
		\rdTo[rightshortfall=2em,crab-]_\cong
	\\
	 \rnode{b}{\int^Y\C(B,Y)\times\C(Y\tn C,Z)}
		&\rTo_{\int^Y\!\!\C(B,Y)\times\curry_{Y,C,Z}}
		&\int^Y\C(B,Y)\times\C(Y,C\lolli Z) && J(B\lolli(C\lolli Z))
	\\
	 \dTo<\cong && \dTo<\cong &\ldTo_{e_{B,C\lolli Z}}
	\\
	 \C(B\tn C,Z) &\lTo_{\curry^{-1}_{B,C,Z}}& \C(B, C\lolli Z)
	\nccurve[angle=180]{->}ab\Aput{\scriptstyle \int^Y\!\!\lambda_{B,Y}\times\C(Y\tn C,Z)}
\end{diagram}
\end{array}\]
\caption{Diagram used in the proof of Lemma~\ref{lemma-2}}\label{fig-1}
\end{sidewaysfigure}
Consider the diagram in Fig.~\ref{fig-1}.
The left-hand cell commutes by the relationship between $e$ and $\lambda$.
The remaining cells commute by naturality,
or functoriality of the coend.
The left edge of this diagram is equal to the left and lower edge
of \pref{bigone}. Therefore condition~\pref{bigone} is
equivalent to the condition in the statement.
\end{proof}
We are still working in an \SECC as originally defined, so we have
a $\curry$ isomorphism. We construct a natural isomorphism $\psi'$ as
follows.
\begin{definition}
	Given a natural isomorphism
	\[
		\curry_{A,B,C}: \C(A\tn B,C)\to\C(A, B\lolli C)
	\]
	we define the natural isomorphism
	\[
		\psi'_{A,B,C}: (A\tn B)\lolli C\enskip\longrightarrow\enskip A\lolli(B\lolli C)
	\]
	to be the unique such natural transformation for which
	\begin{diagram}[h=1.5em]\dlabel{def-psi}
		\C(A,(X\tn Y)\lolli Z) &\rTo^{\C(A,\psi')}& \C(A,X\lolli(Y\lolli Z))\\
		\dTo<{\curry_{A,X\tn Y,Z}^{-1}}\\
		\C(A\tn(X\tn Y),Z)&\dnum&\uTo>{\curry_{A,X,Y\lolli Z}}\\
		\dTo<{\C(\alpha_{A,X,Y},Z)}\\
		\C((A\tn X)\tn Y, Z)&\rTo_{\curry_{A\tn X,Y,Z}}&\C(A\tn X,Y\lolli Z)
	\end{diagram}
	commutes for all $A$,$X$,$Y$,$Z\in\C$.
	(Uniqueness is a consequence of the Yoneda lemma.)
\end{definition}
Using this, we can recast condition~\pref{bigone} very simply.
\begin{lemma}\label{lemma-3}
	In an \SECC, condition~\pref{bigone} holds iff
	\begin{diagram}\dlabel{psi}
		J((A\tn B)\lolli C) &\rTo^{J(\psi'_{A,B,C})}&J(A\lolli(B\lolli C))\\
		\dTo<{e_{A\tn B,C}}&\dnum&\dTo>{e_{A,B\lolli C}}\\
		\C(A\tn B,C)&\rTo_{\curry_{A,B,C}}	&\C(A,B\lolli C)
	\end{diagram}
	commutes for all $A$,$B$,$C\in\C$.
\end{lemma}
\begin{proof}
	By Lemma~\ref{lemma-2} we know that \pref{bigone} is equivalent
	to \pref{diag-e}. Now we have
	\begin{diagram}[h=1.5em,labelstyle=\scriptstyle]
		\JX\C(X\tn(B\tn C),Z) &\rTo^{\JX\C(\alpha_{X,B,C}, Z)}	& \JX\C((X\tn B)\tn C,Z)\\
			&						& \dTo>{\JX\curry_{X\tn B,C,Z}}\\
		\dTo<{\JX\curry_{X,B\tn C,Z}}&\mbox{\pref{def-psi}}& \JX\C(X\tn B,C\lolli Z)\\
			&						& \dTo>{\JX\curry_{X,B,C\lolli Z}}\\
		\JX\C(X,(B\tn C)\lolli Z)&\rTo^{\JX\C(X,\psi'_{B,C,Z})}&\JX\C(X,B\lolli(C\lolli Z))\\
		\dTo<\cong&		\natural& \dTo>\cong\\
		J((B\tn C)\lolli Z) &\rTo_{J(\psi'_{B,C,Z})}& J(B\lolli(C\lolli Z))\\
		\dTo<{e_{B\tn C,Z}}&		& \dTo>{e_{B, C\lolli Z}}\\
		\C(B\tn C,Z)&\rTo_{\curry_{B,C,Z}}& \C(B,C\lolli Z)
	\end{diagram}
	The upper two regions commute for the reasons marked, and all the
	arrows are invertible, therefore the outside \pref{diag-e}
	commutes iff the lower cell \pref{psi} does. 
\end{proof}
\begin{lemma}\label{lemma-4}
	In every \SECC, diagram \pref{psicoh} commutes
	for $\psi'$ in place of $\psi$.
\end{lemma}
\begin{sidewaysfigure}
	\[\begin{array}{@{\hskip8.5em}l}
	\begin{diagram}[labelstyle=\scriptstyle]
		 \rnode{tl}{\C(X, (A\tn(B\tn C))\lolli D)}
			&\rTo^{\C(X,\alpha_{A,B,C}\lolli D)}
			&\rnode{tr}{\C(X,((A\tn B)\tn C)\lolli D))}
		\\
		 \uTo<{\curry_{X,A\tn(B\tn C),D}}&\natural
			&\uTo>{\curry_{X, (A\tn B)\tn C,D}}
		\\
		 \C(X\tn(A\tn(B\tn C)),D)
		 	&\rTo^{\C(X\tn\alpha_{A,B,C},D)}
			&\C(X\tn((A\tn B)\tn C),D)
		\\
		 \dTo<{\C(\alpha_{X,A,B\tn C}, D)}&
		 	\raise-1em\hbox{\pref{pentagon}}
		 	&\dTo>{\C(\alpha_{X,A\tn B,C}, D)}
			&\hbox{\pref{def-psi}}
		\\
		 \C((X\tn A)\tn(B\tn C), D)&
		 	&\C((X\tn(A\tn B))\tn C,D)
		\\
		 &\rdTo(1,2)^{\C(\alpha_{X\tn A,B,C},D)}
		 	\ldTo(1,2)^{\C(\alpha_{X,A,B}\tn C,D)}
		 	&\dTo>{\curry_{X\tn(A\tn B),C,D}}
		\\
		 &\hbox to 10em{\hss$\C(((X\tn A)\tn B)\tn C,D)$}
		 	\hbox to 2em{\hss$\natural$\enskip}
		 	&\C(X\tn(A\tn B), C\lolli D)
			&\pile{
				\rTo^{\curry_{X,A\tn B,C\lolli D}}\\
				\lTo_{\curry_{X,A\tn B,C\lolli D}^{-1}}}
			&\rnode{r}{\C(X, (A\tn B)\lolli(C\lolli D))}
		\\
		 \dTo>{\curry_{X\tn A, B\tn C,D}}&
		 	&\rdTo(1,2)_{\curry_{(X\tn A)\tn B,C,D}}
		 	\dTo>{\C(\alpha_{X,A,B}, C\lolli D)}
		\\
		 &&\C((X\tn A)\tn B, C\lolli D)
		\\
		 &\raise 2em\hbox{\pref{def-psi}}&\dTo>{\curry_{X\tn A,B,C\lolli D}}
		 &\hbox{\pref{def-psi}}
		\\
		 \C(X\tn A, (B\tn C)\lolli D)
		 	&\rTo_{\C(X\tn A,\psi'_{B,C,D})}
			&\C(X\tn A, B\lolli(C\lolli D))
		\\
		 \dTo<{\curry_{X,A,(B\tn C)\lolli D}}&
		 	\natural
		 	&\dTo>{\curry_{X,A,B\lolli(C\lolli D)}}
		\\
		 \rnode{bl}{\C(X,A\lolli((B\tn C)\lolli D))}
		 	&\rTo_{\C(X,A\lolli\psi'_{B,C,D})}
			&\rnode{br}{\C(X,A\lolli(B\lolli(C\lolli D)))}
		\psset{linearc=.5}
		\ncangle[angle=180,arm=30pt]{->}{tl}{bl}
			\Bput{\scriptstyle \C(X,\psi'_{A,B\tn C,D})}
			\Aput{\hskip 3em\hbox{\pref{def-psi}}}
		\nccurve[angleA=0,angleB=90]{->}{tr}r
			\Aput{\scriptstyle \C(X,\psi'_{A\tn B,C,D})}
		\nccurve[angleA=-90,angleB=0]{->}r{br}
			\Aput{\scriptstyle \C(X,\psi'_{A,B,C\lolli D})}
	\end{diagram}
	\end{array}\]
	\caption{Diagram used in the proof of Lemma~\ref{lemma-4}}
	\label{fig-2}
\end{sidewaysfigure}
\begin{proof}
	Consider the diagram in Fig.~\ref{fig-2}. All the regions
	commute for the reasons marked, thus the outside commutes.
	By Yoneda, it follows that \pref{psicoh} commutes as required.
\end{proof}
\begin{lemma}
	Suppose we have a natural isomorphism $\psi$ with
	components
	\[
		\psi_{A,B,C}: (A\tn B)\lolli C \enskip\longrightarrow\enskip A\lolli(B\lolli C)
	\]
	such that
	\begin{diagram}\dlabel{psi-lemma}
		J((A\tn B)\lolli C) &\rTo^{J(\psi_{A,B,C})}&J(A\lolli(B\lolli C))\\
		\dTo<{e_{A\tn B,C}}&\dnum&\dTo>{e_{A,B\lolli C}}\\
		\C(A\tn B,C)&\rTo_{\curry_{A,B,C}}	&\C(A,B\lolli C)
	\end{diagram}
	commutes for all $A$,$B$,$C\in\C$.
	If $\psi$ satisfies condition~\pref{psicoh} then $\psi=\psi'$.
\end{lemma}
\begin{proof}
	Suppose $\psi$ satisfies condition~\pref{psicoh}.
	Then we have
	\begin{mspill}
	\begin{diagram}[h=1.5em,w=2em,labelstyle=\scriptstyle]
		\C(A\tn(B\tn C), D)&&&\rTo^{\C(\alpha_{A,B,C},D)}&&&\C((A\tn B)\tn C,D)\\
		&\rdTo[snake=1em]^{e^{-1}}&&\raise3pt\hbox{$\natural$}
			&&\ruTo[crab+,leftshortfall=-.7em,rightshortfall=1.5em,snake=-1em]^{e}
			\ldTo[crab-,leftshortfall=2em,rightshortfall=-.7em,snake=1em]_{e^{-1}}\\
		&& J(\bigl(A\tn(B\tn C)\bigr)\lolli D) &\rTo^{J(\alpha\lolli D)}
			&J(\bigl((A\tn B)\tn C\bigr)\lolli D)
			&&\dTo>{\curry_{A\tn B,C,D}}\\
		&& &&\dTo>{J(\psi_{A\tn B,C,D})}&\raise1em\hbox to 0pt{\quad\pref{psi-lemma}\hss}\\
		\dTo<{\curry_{A,B\tn C,D}}>{\mbox{\quad\enskip\pref{psi-lemma}}}&&
				\dTo<{J(\psi_{A,B\tn C,D})}&\hbox{\pref{psicoh}}&J((A\tn B)\lolli(C\lolli D))
				&\pile{\rTo\\\lTo}&\C(A\tn B,C\lolli D)\\
		&& &&\dTo>{J(\psi_{A,B,C\lolli D})}&\raise-1em\hbox to 0pt{\quad\pref{psi-lemma}\hss}\\
		&& J(A\lolli\bigl((B\tn C)\lolli D\bigr))&\rTo_{J(A\lolli\psi)}
			&J(A\lolli\bigl(B\lolli(C\lolli D)\bigr))
			&&\dTo>{\curry_{A,B,C\lolli D}}\\
		&\ldTo[crab+,leftshortfall=-.7em,rightshortfall=1.5em,snake=-1em]^{e}
			\ruTo[crab-,rightshortfall=-.7em,leftshortfall=1.5em,snake=2em]_{e^{-1}}
			&&\raise-3pt\hbox{$\natural$}&&
			\luTo[crab+,rightshortfall=-.7em,leftshortfall=1.5em,snake=1em]^{e^{-1}}
			\rdTo[crab-,leftshortfall=-.7em,rightshortfall=2.5em,snake=-2.5em]_{e}
		\\
		\rnode{bl}{\C(A, (B\tn C)\lolli D)}&&
			&\rTo_{\C(A,\psi_{B,C,D})}&&
			&\rnode{br}{\C(A,B\lolli(C\lolli D))}\\
		\ncarc[arcangle=-30]{->}{bl}{br}\Bput{\scriptstyle \C(A,\psi'_{B,C,D})}
	\end{diagram}\end{mspill}\vskip 5em\noindent
	The marked cells commute for the reasons indicated,
	and the outer edge by assumption.
	Since all arrows are invertible, it follows that the
	lower cell also commutes; hence $\psi=\psi'$ by Yoneda.
\end{proof}
In other words, if we are given a $\psi$ satisfying \pref{psicoh},
we can use \pref{psi-lemma} to construct a $\curry$
isomorphism such that \pref{diag-e} holds.
(We have already proved that, given a $\curry$ isomorphism
satisfying \pref{diag-e}, we can use \pref{def-psi} to construct
a natural isomorphism $\psi'$ such that \pref{psicoh} holds.)
This completes the proof of Prop.~\ref{prop-1}.

\subsection{Tensor of Elements}\label{s-tensel}
By Yoneda's lemma, an element of $JA$ corresponds to a
natural transformation $\C(A,-)\To J$. If we have
elements $a\in JA$ and $b\in JB$ then we may define a
natural transformation
\[\begin{array}{r@{}c@{}l}
	\C(A\tn B, -) &\rTo^\cong& \C(A, B\lolli -)\\
		&\rTo^a& J(B\lolli -)\\
		&\rTo^\cong&  \C(B, -)\\
		&\rTo^b&  J,
\end{array}\]
corresponding to an element of $J(A\tn B)$. We denote
this element $a\tn b$.
It is easy to check that this operation defines a natural
transformation with components $m_{A,B}: JA\times JB\to J(A\tn B)$.

\begin{propn}
This natural transformation agrees with the associativity,
in the sense that the diagram
\begin{diagram}[h=1.5em]
	\C(A\tn (B\tn C), -)\\
	&\rdTo[snake=1.3em]^{a\tn (b\tn c)}\\
	\dTo<{\C(\alfa ABC, -)} && J\\
	&\ruTo[snake=1.3em]_{(a\tn b)\tn c}\\
	\C((A\tn B)\tn C, -)
\end{diagram}
commutes for every $A$,$B$,$C\in\C$ with
$a\in JA$, $b\in JB$ and $c\in JC$.
\end{propn}
\begin{proof}
	Consider the diagram
	\begin{diagram}
	\C(A\tn(B\tn C),Z)&\rTo^\cong&\C(A,(B\tn C)\lolli Z)&\rTo^{a}&J\bigl((B\tn C)\lolli Z)
		&\rTo^\cong&\C(B\tn C,Z)\\
	&& &\rdTo(2,4)_{\C(A,\psi)}& &\rdTo(2,4)_{J(\psi)}\hbox to0em{\pref{psi}}& \dTo>\cong\\
	\dTo<{\C(\alpha,Z)} &\pref{def-psi}& && \natural && \C(B,C\lolli Z)\\
	&& && && \uTo>\cong\\
	\C((A\tn B)\tn C,Z)&\rTo_\cong&\C(A\tn B,C\lolli Z)&\rTo_\cong&\C\bigl(A, B\lolli(C\lolli Z)\bigr)
		&\rTo_{a}&J\bigl(B\lolli(C\lolli Z)\bigr)\\
	\end{diagram}
	Since the internal cells commute, so does the outside.
	Now observe that, by definition, the natural transformation
	$a\tn(b\tn c)$ is equal to the upper path followed by the
	composite
	\[
		\C(B,C\lolli Z) \rTo^b J(C\lolli Z) \cong \C(C,Z) \rTo^c JZ,
	\]
	while $(a\tn b)\tn c$ is equal to the lower path followed
	by this composite. Thus the claim follows.
\end{proof}



\section{The Star-autonomous Case}\label{s-sstac}
There is a general notion of promonoidal \staraut category
\citep[][\S7]{DSQuantum}. A symmetric promonoidal \staraut category
is a symmetric promonoidal category $\C$ equipped with a full and
faithful functor $-\perp: \C\to\C\op$ and a natural
isomorphism
\[
	P(A, B, C\perp) \cong P(A, C, B\perp).
\]
This specialises in the obvious way:
\begin{definition}
	A \sstac is a \ssmc $\C$ with a
	full and faithful functor $-\perp:\C\to\C\op$
	and a natural isomorphism
	\begin{equation}\label{eq-sstac}
		\C(A\tn B, C\perp) \cong \C(A\tn C, B\perp).
	\end{equation}
\end{definition}
Note that, since $-\perp$ is full and faithful, there is a natural
isomorphism $\C(A,B) \cong \C\op(A\perp,B\perp) = \C(B\perp, A\perp)$.

\begin{lemma}\label{lemma-perpp}
	There is a natural isomorphism $B\cong B\perpp$.
\end{lemma}
\begin{proof}
	There is a sequence of natural isomorphisms
	\[\begin{array}{r@{{}\protect\cong{}}lp{5cm}}
		\C(A,B) & \C(B\perp, A\perp)	&$\perp$ is full and faithful\\
		& \int^XJX\times\C(X\tn B\perp, A\perp) &using $\lambda^{-1}$\\
		& \int^XJX\times\C(X\tn A, B\perpp)	    &by \pref{eq-sstac}\\
		& \C(A, B\perpp)				&using $\lambda$.
	\end{array}\]
	Therefore, by Yoneda's lemma, it follows that
	$B$ is naturally isomorphic to $B\perpp$, as required.
\end{proof}
\begin{propn}
	Every \sstac is a \SECC, with $A\lolli B$ defined
	as $(A\tn B\perp)\perp$.
\end{propn}
\begin{proof}
	Clearly $\lolli$ is a functor of the correct type,
	so it remains only to establish the existence of
	a natural isomorphism $\C(A\tn B,C)\cong\C(A,B\lolli C)$.
	We have the following sequence of isomorphisms:
	\[\begin{array}{r@{\;}c@{\;}lp{5cm}}
		\C(A\tn B, C) &\cong& \C(A\tn B, C\perpp) & by Lemma~\ref{lemma-perpp}\\ 
		&\cong& \C(B\tn A, C\perpp) & by symmetry\\
		&\cong& \C(B\tn C\perp, A\perp) & by \pref{eq-sstac}\\
		&\cong& \C(A\perpp, (B\tn C\perp)\perp) &since $\perp$ is full and faithful\\
		&\cong& \C(A, (B\tn C\perp)\perp) & by Lemma~\ref{lemma-perpp}\\
		&=&     \C(A, B\lolli C) &by definition of $\lolli$.
	\end{array}\]
	\vskip-2em
\end{proof}

\section{Related Work}\label{s-related}
\subsection{The Lamarche-Stra{\ss}burger Definition}\label{s-ls}

Not long after starting work on this we came
across the draft of \citet{LSFreeBool}. The authors define what they
call \emph{(unitless) autonomous categories}, motivated (like us) by
the desire to model unitless fragments of MLL.


Using our notation, an `autonomous category' in the sense of
\citet{LSFreeBool} consists of:
\begin{itemize}
        \item a category $\C$,
        \item functors $\tn:\C\times\C\to\C$ and $\lolli:\C\op\times\C\to\C$
                with a natural isomorphism
                \(
                \curry_{A,B,C}: \C(A\tn B,C) \rTo^\cong \C(A, B\lolli C),
                \)
        \item natural isomorphisms $\alpha$ and $\sigma$ such that
                $\sigma_{A,B} = \sigma_{B,A}^{-1}$,
                satisfying conditions \pref{pentagon} and \pref{hexagon},
        \item a functor $J:\C\to\Set$, with a natural isomorphism
                \[
                e_{A,B}: J(A\lolli B) \rTo^\cong \C(A,B),
                \]
                subject to the condition obtained by translating \pref{II}
                into the promonoidal setting.
\end{itemize}
(For the original definition in the draft \citep{LSFreeBool}, in terms
of \emph{virtual objects}, see Appendix~\ref{app-ls}.)
It is clear that every \SECC is a Lamarche-Stra\ss burger category;
however the converse appears to be false. In fact, certain properties are
desired of these categories which do not appear to be derivable from
the given conditions. In the case where the functor $J$ is representable,
the axioms do not seem to imply that the category is symmetric monoidal
closed.
Also it is claimed (p.~3 of op. cit.)
that the canonical natural transformation $JA\times JB \to J(A\tn B)$
`agrees well with associativity'.\footnote{
	Our definition does have this property, as
	shown in \S\ref{s-tensel}.
}
In fact, one of the motivations
behind producing this preprint is to point out that, at least, we can
provide a solution to the problem of finding a definition with the
desired properties. In correspondence Lamarche and Stra\ss
burger have indicated that they might change this definition in the
final version of their paper.

For the moment, we should like to point out that their current
definition of `autonomous functor' (their Def.~2.1.4) also does not
quite do what one might expect.
On the one hand it does not demand
(nor imply) that the functor preserves the $\curry$ isomorphism, on
the other it demands a natural \emph{isomorphism} $J\cong JF$. (This
latter condition is not satisfied, for example, by the unique functor
$\Set\to\mathbf1$.)

\subsection{The Do{\v s}en-Petri{\'c} Definition}\label{s-dp}
\citet{ProofNetCats} define what they call a \defn{proof-net category}
to be a unitless linearly distributive category in which each object
has a dual in a suitable sense. This is a reasonable approach, and
we conjecture that the resulting definition is equivalent to ours.

Early in our exploration of candidates for \sstac,
we considered essentially the same definition.
Ultimately we chose the approach more analogous to the
standard progression from monoidal category to \stac,
via symmetry and closure: it ties directly into the pioneering work
of \citet{EKClosed}, with the \pref{psicoh} diagram, and also
Day's promonoidal categories \citeyearpar{DayPro}.
Along the way it gives us a resonable notion of model for
the intuitionistic fragment of \mll.

\subsection{Eilenberg and Kelly}\label{s-ek}
The pioneering work of \citet{EKClosed} -- in which closed
categories are defined for the first time -- also deserves to be
mentioned here. The authors define `symmetric monoidal closed
category' using a large number of axioms with a great deal
of redundancy. (Our diagram \pref{psicoh} is one of them.)
It seems reasonable to conjecture that, if one were to delete
the unit and the axioms involving it from this original
definition, the structures satisfying the remaining axioms
would be just the \SECCs.

\section{Ongoing Work}
There are general definitions of \emph{lax promonoidal functor}
\citep{DayMonoidalMonads}, \emph{strong promonoidal functor}
\citep{PromonKan}, and \emph{promonoidal \staraut{} functor}
\citep{DSQuantum}.
These definitions need to be specialised to \SECCs and \sstacs in a
sensible way.

We are also working on refining the definition of \sstac, to
give an elementary description that does not rely on the
definition of \SECC.
We conjecture that our definition is equivalent to
\citeauthor{ProofNetCats}'s proof-net categories -- see \S\ref{s-dp}.

As mentioned in the introduction, we hope to show that the proof-net
category of \citet{MallNets} is the free \sstac with finite products
(free in a 2-categorical sense\footnote{Since MALL formulas are implicitly
quotiented by de Morgan duality.}).

\appendix
\section{Axioms for Monoidal Categories}\label{app-ax}
\renewcommand\theequation{\thesection.\arabic{equation}}
\setcounter{equation}{0}
This appendix describes some axioms for a
symmetric monoidal category, and shows that this axiomatisation
is equivalent to the usual one.\footnote{
	There is every chance that this axiomatisation is somewhere
	in the literature. I don't know a reference for it though.
}
(Of course we are really interested in promonoidal categories:
the arguments here can readily be transferred to the more
general setting.)

A monoidal category is a category $\C$ equipped with a functor
\[
	\tn:\C\times\C \to\C
\]
and an object $I\in\C$, together with natural isomorphisms
$\alpha$, $\lambda$ and $\rho$ having components
\[\begin{array}{r@{{}\to{}}l}
	\alfa ABC: (A\tn B)\tn C & A\tn(B\tn C)\\
	\lambda_A: I\tn A & A\\
	\rho_A: A\tn I & A
\end{array}\]
such that the following diagrams commute for all $A$,$B$,$C$,$D\in\C$:

\begin{mspill}\begin{diagram}
  \bigl((A\tensor B)\tensor C\bigl)\tensor D
  &\rTo^{\alpha_{A\tn B,C,D}}&(A\tensor B)\tensor (C\tensor D)
  &\rTo^{\alpha_{A,B,C\tn D}} & A\tensor \bigl(B\tensor (C\tensor D)\bigr)
  \\
  &\rdTo[snake=-1em](1,2)<{\alpha_{A,B,C}\tn D} &\dnum[pentagon]
  & \ruTo[snake=1em](1,2)>{A\tn\alpha_{B,C,D}}
  \\
  & \spleft{\bigl(A\tensor(B\tensor C)\bigr)\tensor D}
  & \rTo_{\alpha_{A,B\tn C,D}}
  & \spright{A\tensor\big((B\tensor C)\tensor D\big)}
\end{diagram}\end{mspill}

\begin{diagram}\dlabel{AIC}
	(A\tn I)\tn C &\rTo^{\alpha_{A,I,C}}&A\tn(I\tn C)\\
	&\rdTo[snake=-1ex](1,2)<{\rho_A\tn C}\raise1ex\dnum\ldTo[snake=1ex](1,2)>{A\tn\lambda_C}\\
	&A\tn C
\end{diagram}
These axioms have many interesting consequences. Most importantly,
it follows that the following three diagrams commute for all
$A$,$B$,$C\in\C$:
\[\begin{array}{c@{\qquad}c}
	\begin{diagram}\dlabel{IBC}
		(I\tn B)\tn C &\rTo^{\alfa IBC}&I\tn(B\tn C)\\
		&\rdTo[snake=-1ex](1,2)_{\lambda_B\tn C}\raise1ex\dnum
			\ldTo[snake=1ex](1,2)_{\lambda_{B\tn C}}\\
		&B\tn C
	\end{diagram}
	&
	\begin{diagram}\dlabel{ABI}
		(A\tn B)\tn I &\rTo^{\alfa ABI}&A\tn(B\tn I)\\
		&\rdTo[snake=-1ex](1,2)_{\rho_{A\tn B}}\raise1ex\hbox{(\theequation)}\ldTo[snake=1ex](1,2)_{A\tn\rho_B}\\
		&A\tn B
	\end{diagram}
	\\
	\multicolumn 2c{
	\begin{diagram}\dlabel{II}
		\rnode{II}{I\tn I} &\dnum& \rnode{I}{I}\\
		\nccurve[angleA=60,  angleB=120] {->}{II}{I}\Aput{\lambda_I}
		\nccurve[angleA=-60, angleB=-120]{->}{II}{I}\Bput{\rho_I}
	\end{diagram}}
\end{array}\]
\citet{BTC} give a simple and elegant proof.
There are other possible axiomatisations: for example, conditions
\pref{pentagon}, \pref{IBC} and \pref{II} are also collectively
sufficient.\footnote{This can be proved by the technique of \citet{BTC}.}

A \emph{symmetric} monoidal category is a monoidal category $\C$
equipped with a natural isomorphism $\sigma$ having components
\[
	\sigma_{A,B}: A\tn B \to B\tn A,
\]
such that $\sigma_{A,B} = \sigma_{B,A}^{-1}$ and the following
commutes for all $A$,$B$,$C\in\C$:
\begin{diagram}
  (A\tn B)\tn C &\rTo^{\alfa ABC} & A\tn(B\tn C)
    & \rTo^{\sigma_{A,B\tn C}} & (B\tn C)\tn A \\
  \dTo<{\sigma_{A,B}\tn C} &&\dnum[hexagon]&& \dTo>{\alfa BCA}\\
  (B\tn A)\tn C & \rTo_{\alpha_{B,A,C}} & B\tn(A\tn C)
    & \rTo_{B\tn\sigma_{A,C}} & B\tn(C\tn A)
\end{diagram}
It then follows that
\begin{diagram}[w=4em]\dlabel{IA}
	I\tn A &\rTo^{\sigma_{I,A}}&A\tn I\\
	&\rdTo[snake=-1ex](1,2)_{\lambda_A}
	\raise1ex\dnum
	\ldTo[snake=1ex](1,2)_{\rho_A}\\
	&A
\end{diagram}
commutes for every $A\in\C$; again, see \citet{BTC} for a proof.

Conditions \pref{pentagon} and \pref{hexagon} are indispensable;
the usual definition of symmetric monoidal category requires \pref{AIC}
in addition. However, it is sometimes convenient to eliminate the
$\rho$ isomorphism from the data: that is permissible, since
\pref{IA} shows that $\rho$ may be defined in terms of $\lambda$
and $\sigma$. It turns out that, in this situation, we may require
\pref{IBC} in place of \pref{AIC}. Specifically:
\begin{propn}\label{prop-A1}
	If \pref{IBC}, \pref{hexagon} and \pref{IA} hold then so does \pref{AIC}.
\end{propn}
\begin{proof}
Consider the diagram
\let\nat = \natural%
\begin{diagram}[labelstyle=\scriptstyle]
  &&(I\tn A)\tn C &\rput(0,-.5){\hbox{\pref{IBC}}}\rTo^{\alfa IAC}&I\tn(A\tn C)
 \\
  &\ruTo^{\sigma_{A,I}\tn C}&\dTo>{\lambda_A\tn C}
  &\ldTo_{\lambda_{A\tn C}}&\nat&\rdTo^{I\tn\sigma_{A,C}}
 \\
  (A\tn I)\tn C&\rTo_{\rho_A\tn C}&\rput(-.4,.9){\hbox{\pref{IA}}}A\tn C
  &\rTo_{\sigma_{A,C}}&C\tn A\rput(.5,-.6){\hbox{\pref{IBC}}}
  &\lTo^{\lambda_{C\tn A}}&I\tn(C\tn A)
 \\
  &\rdTo_{\alfa AIC}&\uTo>{A\tn\lambda_C}&\nat&\uTo<{\lambda_C\tn A}&\ruTo_{\alfa ICA}
 \\
  &&A\tn(I\tn C)&\rTo_{\sigma_{A,I\tn C}}&(I\tn C)\tn A
\end{diagram}
The outside is an instance of \pref{hexagon}, and the labelled
regions commute for the reasons marked.
Since all the morphisms are invertible, it follows that the
unlabelled region at lower left commutes. This region is just \pref{AIC}.
\end{proof}
In summary, we may define a symmetric monoidal category to be
a category $\C$ with a functor $\tn$ and a unit object $I$,
together with natural isomorphisms $\alpha$, $\lambda$ and
$\sigma$ such that $\sigma_{A,B} = \sigma_{B,A}^{-1}$ and
diagrams \pref{pentagon}, \pref{hexagon} and \pref{IBC} commute.

\section{The Lamarche-Stra\ss burger Definition}\label{app-ls}
\newcommand\A{\mathbb{A}}
\newcommand\B{\mathbb{B}}
\newcommand\I{\mathbb{I}}
\renewcommand\C{\mathcal{C}}

%

In their draft, \cite{LSFreeBool} give a very interesting discussion
of autonomous categories without units, and the following definition
of \emph{(unitless) autonomous category}, based on the notion of a
\emph{virtual object}.  In Section~\ref{s-ls} we presented the definition
in our own notation (i.e., promonoidal style); for the sake of
completeness, here is (a condensed presentation of) the original
definition of \cite{LSFreeBool}, in terms of virtual objects.

A category $\C$ \emph{has tensors} if it is
equipped with a bifunctor $-\tn-:\C\times\C\to\C$ with the usual
associativity and symmetry isomorphisms $$\begin{array}{rl}
{}\textsf{assoc}_{A,B,C}:& A \tn(B \tn C) \to (A \tn B) \tn C\\
{}\textsf{twist}_{A,B} : & A \tn B \to B \tn A
\end{array}$$
obeying the usual associated coherence laws.
When it exists, write $-\lolli-:\C\op\times\C\to\C$ for internal hom,
defined by adjointness to tensor as in the usual case with units.
Write $h_X=\C(-,X)$ and $h^X=\C(X,-)$ for the hom functors, and
$H^X=X\lolli-:\C\to \C$.
Writing a functor $\C\to\Set$ as $h^\A$ for a symbol $\A$ allows one
to write an element $s$ of the set $h^\A(s)$, corresponding (by
Yoneda's Lemma) to a natural transformation $h^X\to h^\A$, as
\[
	\A \rTo.^s X
\]
When the symbol $\A$ is an object of $\C$, the functor $h^\A$ is
representable, in the usual sense; when $\A$ is not an object of $\C$,
it is a \emph{virtual object}.
In general a dotted arrow will mean at least one of the source or
target is virtual, and should be interpreted as a reverse-direction
natural transformation between the corresponding functors. For
example, given $f : X \to Y$ and $t = (h^{\A}f)(s)$, one can draw the
`commutative diagram'
\begin{diagram}[h=1em,w=1.5em]
	&&\A\\
	&\ldTo.^s && \rdTo.^t\\
	X && \rTo_f && Y\;\;,
\end{diagram}
justifying the notation $t = f \circ s$, or simply $t = fs$.  Define
$\A\tn X$ in the obvious way, i.e., $h^{\A\tn X}=h^\A H^X$.  This
construction is natural in both variables: given $s:\A \rTo. \B$
(between virtual objects) and a morphism $f:X\to Y$, there is an
obvious \(s\tn f:\A\tn X \rTo. \B\tn Y\).

\begin{definition}
A category $\C$ with tensors is an \emph{autonomous category} if it
has an internal hom $\lolli$ and a functor $h^\I$ with a natural isomorphism
\[h^\I(X\lolli Y)\;\cong\;\C(X,Y)\]
such that the following diagram (of mostly virtual arrows) commutes:
\begin{diagram}[h=1.5em]
	&&X &\rTo.^\cong&\I\tn X&\rTo.^{t\tn X}& Y\tn X\\
	&\ruTo.^s\\
	\I && && && \dTo>\cong\\
	&\rdTo._t\\
	&&Y&\rTo._\cong&\I\tn Y&\rTo._{s\tn Y}&X\tn Y
\end{diagram}
\end{definition}
This diagram corresponds to diagram \pref{II} of the previous
appendix, translated into promonoidal style.

\bibliography{cs}

\begin{thebibliography}{15}
\expandafter\ifx\csname natexlab\endcsname\relax\def\natexlab#1{#1}\fi
\expandafter\ifx\csname url\endcsname\relax
  \def\url#1{{\tt #1}}\fi

\bibitem[Barr(1979)]{BarrStac}
Michael Barr.
\newblock {\em {$\star$}-Autonomous Categories}, volume 752 of {\em Lecture
  Notes in Mathematics}.
\newblock Springer-Verlag, 1979.

\bibitem[C{\`a}ccamo et~al.(2002)C{\`a}ccamo, Hyland, and Winskel]{CatNotes}
Mario C{\`a}ccamo, Martin Hyland, and Glynn Winskel.
\newblock Lecture notes in category theory.
\newblock URL \url{http://www.brics.dk/~mcaccamo/catnotes.ps.gz}.
\newblock BRICS LS-02-catnotes (draft), 2002.

\bibitem[Day(1970)]{DayPro}
Brian Day.
\newblock On closed categories of functors.
\newblock In Saunders Mac~Lane, editor, {\em Reports of the Midwest Category
  Seminar}, volume 137 of {\em Lecture Notes in Mathematics}, pages 1--38,
  Berlin--New York, 1970. Springer-Verlag.

\bibitem[Day(1977)]{DayMonoidalMonads}
Brian Day.
\newblock Note on monoidal monads.
\newblock {\em Journal of the Australian Mathematical Society}, 23:\penalty0
  292--311, 1977.

\bibitem[Day and Street(1995)]{PromonKan}
Brian Day and Ross Street.
\newblock Kan extensions along promonoidal functors.
\newblock {\em Theory and Applications of Categories}, 1\penalty0 (4):\penalty0
  72--77, 1995.

\bibitem[Day and Street(1997)]{MonBicat}
Brian Day and Ross Street.
\newblock Monoidal bicategories and {H}opf algebroids.
\newblock {\em Advances in Mathematics}, 129:\penalty0 99--157, 1997.

\bibitem[Day and Street(2004)]{DSQuantum}
Brian Day and Ross Street.
\newblock Quantum categories, star autonomy, and quantum groupoids.
\newblock In {\em Galois Theory, {H}opf Algebras, and Semiabelian Categories},
  volume~43 of {\em Fields Institute Communications}, pages 187--226. American
  Mathematical Society, 2004.
\newblock {\tt arXiv:math.CT/0301209}.

\bibitem[Do{\v s}en and Petri{\'c}(2005)]{ProofNetCats}
Kosta Do{\v s}en and Zoran Petri{\'c}.
\newblock Coherence of proof-net categories.
\newblock {\tt arXiv:math.CT/0503301}, March 2005.

\bibitem[Eilenberg and Kelly(1965)]{EKClosed}
Samuel Eilenberg and Max Kelly.
\newblock Closed categories.
\newblock In {\em Conference on Categorical Algebra, La Jolla}, pages 421--562.
  Springer-Verlag, 1965.

\bibitem[Eilenberg and Kelly(1966)]{EK66}
Samuel Eilenberg and Max Kelly.
\newblock A generalization of the functorial calculus.
\newblock {\em J.\ Algebra}, 3:\penalty0 366--375, 1966.

\bibitem[Girard(1987)]{Girard87}
Jean-Yves Girard.
\newblock Linear logic.
\newblock {\em Theoretical Computer Science}, 50:\penalty0 1--102, 1987.

\bibitem[Hughes and van Glabbeek(2005)]{MallNets}
Dominic Hughes and Rob van Glabbeek.
\newblock Proof nets for unit-free multiplicative-additive linear logic.
\newblock {\em ACM Transactions on Computational Logic}, 2005.
\newblock To appear. Invited submission November 2003, revised January 2005.

\bibitem[Joyal and Street(1993)]{BTC}
Andr{\'e} Joyal and Ross Street.
\newblock Braided tensor categories.
\newblock {\em Advances in Mathematics}, 102\penalty0 (1):\penalty0 20--78,
  1993.

\bibitem[Kelly and Mac~Lane(1971)]{KM71}
Max Kelly and Saunders Mac~Lane.
\newblock Coherence in closed categories.
\newblock {\em J.\ Pure and Applied Algebra}, 1:\penalty0 97--140, 1971.

\bibitem[Lamarche and Stra{\ss}burger(2005)]{LSFreeBool}
Fran{\c c}ois Lamarche and Lutz Stra{\ss}burger.
\newblock Constructing free boolean categories.
\newblock Submission for \emph{Logic in Computer Science}, January 2005.

\end{thebibliography}
\end{document}